\documentclass[A4paper,11pt,english]{article}
\usepackage{babel}
\usepackage{graphicx}
\usepackage{amsmath}
\usepackage[utf8]{inputenc}
\usepackage{amsfonts}
\usepackage{amssymb}
\usepackage{amsmath}
\usepackage{amsthm}
\usepackage{xcolor}
\usepackage{array}
\usepackage{hyperref}
\title{\bf Max-product Kantorovich sampling operators: quantitative estimates in functional spaces}
\author{ {\bf Lorenzo Boccali} \hskip1cm {\bf Danilo Costarelli} \hskip1cm {\bf Gianluca Vinti} \\ 
	Department of Mathematics and Computer Science \\
	University of Perugia\\
	1, Via Vanvitelli, 06123 Perugia, Italy \\ 
	{\small {\tt lorenzo.boccali@unipg.it}} - {\small {\tt danilo.costarelli@unipg.it}} - {\small {\tt gianluca.vinti@unipg.it}} }
\date{}

\newtheorem{lemma}{Lemma}[section]
\newtheorem{teorema}{Theorem}[section]
\newtheorem{cor}{Corollary}[section]
\theoremstyle{definition}
\newtheorem{remark}{Remark}[section]
\newcommand{\N}{\mathbb{N}}
\newcommand{\R}{\mathbb{R}}
\newcommand{\Z}{\mathbb{Z}}
\newcommand{\la}{\lambda}
\newcommand{\fhi}{\varphi}

\newcommand{\Om}{\Omega}
\newcommand{\assolutol}{\bigl\lvert}
\newcommand{\assolutor}{\bigr\rvert}
\newcommand{\norma}{\|}

\newcommand{\I}{I^{\varphi}}
\newcommand{\op}{K_{n}^{\chi}}
\newcommand{\pop}{P_{n}^{\chi}}
\newcommand{\mom}{m_{0}(\chi)}

\newcommand{\integrale}{\int_{\Omega}}
\newcommand{\integralek}{\int_{k/n}^{(k+1)/n}}
\newcommand{\integraleab}{\int_{a}^{b}}
\newcommand{\integraler}{\int_{\R}}

\newcommand{\Vz}{\bigvee_{k \in \mathbb{Z}}}
\newcommand{\Vj}{\bigvee_{k \in \mathcal{J}_{n}}}

\begin{document}
	\maketitle
	\begin{abstract}
\noindent In this paper, we study the order of approximation for max-product Kantorovich sampling operators based upon generalized kernels in the setting of Orlicz spaces. We establish a quantitative estimate for the considered family of sampling-type operators using the Orlicz-type modulus of smoothness, which involves the modular functional of the space. From this result, it is possible to obtain the qualitative order of convergence when functions belonging to suitable Lipschitz classes are considered. On the other hand, in the compact case, we exploit a suitable definition of K-functional in Orlicz spaces in order to provide an upper bound for the approximation error of the involved operators. The treatment in the general framework of Orlicz spaces allows one to obtain a unifying theory on the rate of convergence, as the proved results can be deduced for a wide range of functional spaces, such as $L^{p}$-spaces, interpolation spaces and exponential spaces.    
\end{abstract}
\medskip\noindent
{\small {\bf AMS subject classification:} 41A25, 41A35  \newline
{\small {\bf Key Words:} Max-product Kantorovich sampling operators; quantitative estimates; Orlicz spaces; modulus of smoothness; K-functional; Lipschitz classes.
\section{Introduction}
Sampling theory, with its applications to modern disciplines such as information theory and communication engineering (see, e.g., \cite{Benedetto,Marvasti}), has seen a continuous increase of interest in the last years. One of the crucial facts in that theory has been the introduction of a generalized version of the classical cardinal series, thanks to P. L. Butzer and his school in Aachen around the $1980$s. The main goal of Butzer's research was to develop a theory that would overcome the limitations, from an application point of view, of the celebrated Whittaker-Kotel'nikov-Shannon (WKS) sampling theorem (see, e.g., \cite{Butzer,Higgins,Higgins2}). This need was due to the fact that real-world signals do not generally achieve the high degree of regularity required by the assumptions of the above theorem. For this reason, many generalizations of the WKS theorem have been introduced and studied from then until now, in order to consider the approximation of both continuous and not-necessarily continuous functions (see, e.g., \cite{Ries,Butzer2,Butzer3,Kivinukk,Kivinukk2}). In particular, in approximation processes concerning the latter kind of signals, the Kantorovich approach, introduced in the general context of Orlicz spaces by Bardaro et al. \cite{Bardaro2}, has been shown to be very successful. More precisely, Kantorovich-type operators are motivated by the idea that, in concrete cases, the behavior of a signal can be better caught if it is sampled in a whole neighborhood of a given point rather than at a single sampling node. This fact is due to the presence of what, in Signal Analysis, it is called time-jitter error}. From a mathematical point of view, the Kantorovich approach can be translated into replacing the sample values $f(k/w)$ in the definition of a family of discrete sampling operators by mean values of the form  $w \int_{k/w}^{(k+1)/w} f(t) \ dt$, $k \in \Z$, $w >0$, returning a (linear) new version of sampling series that have good approximation properties in $L^{p}$-spaces, with $1 \le p < +\infty$ (see, e.g., \cite{Orlova,Costarelli,CV2}). Such modification, in fact, reduces the influence of the time-jitter error in the reconstruction process. 
\newline More recently, in \cite{Coroianu}, the authors introduced and studied the max-product version of Kantorovich sampling operators based upon generalized kernels acting on the $L^{p}$-spaces and suitable spaces of continuous functions. In Operator Theory, the max-product approach, firstly introduced by Bede, Coroianu and Gal (for a very detailed survey, see the monograph \cite{Bede}) can be described as a generalization process that transforms any family of linear approximation operators, defined by finite sums (or series in the case of infinite terms), into a non-linear one, constructed by computing the maximum (or supremum, in the infinite case) over the same set of indexes (see, e.g., \cite{Coroianu3,Gungor,Holhos,Coroianu2,CSV}). The main advantage of max-product operators is essentially that they have a better order of approximation than their linear counterparts. In recent years, this has led several researchers to focus their attention on this type of non-linear operators, proving strong localization results and shape-preserving properties (see, e.g., \cite{Bede2,Bede3,CG2,CG,CV}). In addition, max-product type operators find interesting applications not only to function approximation, but also in fuzzy sets theory and possibility theory (see, e.g., \cite{Bede4,Gal,Cornejo}). 
\newline In this paper, we consider the family $(\op)_{n \in \N}$ of max-product sampling-type operators defined by: \[\op(f)(x):=\dfrac{\displaystyle \Vz \chi(nx-k) \left[n \integralek f(t) \ dt\right]}{\displaystyle \Vz \chi(nx-k)}, \quad x \in \R, \quad \textnormal{(I)},\]
where $f: \R \rightarrow \R$ is a locally integrable function and $\chi: \R \rightarrow \R$ is a generalized kernel function satisfying suitable assumptions (see, e.g., \cite{CCGV,Coroianu}). For the sake of completeness, we stress that it is not necessarily assume that $\chi$ is a (discrete) approximate identity (see, e.g., \cite{Cruz}), but convergence results for $\op$ can be proved under weak conditions.  Here, the symbol $\bigvee$ denotes the supremum of the corresponding set of real numbers. \newline In \cite{Boccali}, we proved a modular convergence theorem for the above operators in the general setting of Orlicz spaces generated by convex $\fhi$-functions, which include the Lebesgue spaces and other well-known functional spaces as particular cases. The approximation of functions defined on both bounded intervals and on the whole real axis have been considered. \newline For what concerns the rate of convergence, in \cite{Coroianu}, the order of approximation for (I) has been studied in the classical space of uniformly continuous and bounded functions, i.e., $C(\R)$ endowed with the usual sup-norm $\norma \cdot \norma_{\infty}$, and in $L^{p}(\R)$, in both cases using a suitable version of the modulus of continuity of the function to be approximated. In the compact case, a quantitative estimate with respect to the norm $\norma \cdot \norma_{p}$ has been obtained via the well-known Peetre K-functional (\cite{Peetre}).
\newline The purpose of this paper is to investigate the quantitative modular convergence properties of the max-product Kantorovich sampling operators, i.e., we want to obtain an upper bound for the modular of the error of approximation using suitable moduli of smoothness and K-functionals in Orlicz spaces. To do this, in Section \ref{preliminari}, we recall the notion of modulus of smoothness $\omega(f, \delta)_{\fhi}$ in $L^{\fhi}(\R)$ \cite{Bardaro3}, (defined in terms of the modular functional of the space) together with the definition of the considered operators and other auxiliary results. In Section \ref{risultatiprincipali}, we prove quantitative estimates for the approximation of functions $f:\R \rightarrow \R$ in terms of $\omega(f, \delta)_{\fhi}$, and via a suitable Orlicz-type K-functional, if $f$ is defined on a bounded interval. Finally, we recall some special instances of Orlicz spaces for which the above results can be applied.      
\section{Preliminary notions}
\label{preliminari}
We begin this section recalling some basic notions needed to introduce the general setting of Orlicz spaces in which we will work.  \newline In what follows, we will call $\fhi$-function any function $\fhi: \R_{0}^{+} \rightarrow \R_{0}^{+}$ which satisfies the following conditions: \vskip0.1cm $(\fhi1)$ $\fhi$ is a non-decreasing and continuous function;
\vskip0.1cm $(\fhi2)$ $\fhi(0)=0$, $\fhi(u)>0$ if $u>0$ and $\displaystyle \lim_{u \rightarrow +\infty} \fhi(u)=+\infty$. \vskip0.2cm \noindent Now, let $\fhi$ be a fixed $\fhi$-function. Then we can introduce the functional $\I:M(\Om)\rightarrow [0, +\infty]$ defined by: 
\[\I[f]:=\integrale\fhi\bigl(\assolutol f(x) \assolutor\bigr) \ dx,\]
for every $f \in M(\Om)$, where $M(\Om)$ denotes the set of all (Lebesgue-)measurable functions $f: \Om \rightarrow \R$, with $\Om=[a,b]$ or $\Om=\R$. It is well known from the literature (see, e.g., \cite{Musielak,Musielak2,Rao,Bardaro3}) that $\I$ is a modular functional and the corresponding modular space, called the Orlicz space generated by $\fhi$, is defined as follows:
\[L^{\fhi}(\Om):=\bigl\{f \in M(\Om): \I\bigl[\la f\bigr]< +\infty, \ \textnormal{for some} \ \la>0\bigr\}.\] 
The definition of the space $L^{\fhi}(\Om)$ is due to the Polish mathematician W. Orlicz (see \cite{Musielak}) and such theory has been extensively developed by J. Musielak, and several other authors (see, e.g., \cite{Musielak2,Bardaro3}). From now on, we denote by $L_{+}^{\fhi}(\Om)$ the subspace of $L^{\fhi}(\Om)$ of the non-negative functions. 
\newline In order to study approximation results in $L^{\fhi}(\Om)$, a well-known notion of convergence, i.e., the so-called modular convergence, can be introduced using the definition of the modular functional $\I$ (\cite{Musielak}). A net of functions $(f_{w})_{w>0} \subset L^{\fhi}(\Om)$ is said to be modularly convergent to a function $f \in L^{\fhi}(\Om)$ if: 
\begin{equation}
\label{convergenzamodulare}
\I\bigl[\la\bigl(f_{w}-f\bigr)\bigr]=\integrale \fhi\bigl(\la \assolutol f_{w}(x)-f(x) \assolutor\bigr) \ dx \longrightarrow 0, \quad w\rightarrow+\infty,
\end{equation}
for some $\la >0$. The definition given in (\ref{convergenzamodulare}) induces a topology in the space, called modular topology.
\newline Now, in order to establish quantitative estimates for the rate of approximation for a family of max-product sampling-type operators, we need to recall the definition of the Orlicz-type modulus of smoothness. For the sake of simplicity, we consider the case of $\Om=\R$. \newline For any fixed $f \in L^{\fhi}(\R)$, we denote by $\omega(f, \delta)_{\fhi}$ the modulus of smoothness of $f$, defined by:
\[\omega(f,\delta)_{\fhi}:=\sup_{\lvert h \rvert \le \delta}\I[f(\cdot+h)-f(\cdot)],\]
with $\delta >0$. It is well-known (see \cite{Bardaro3}, Theorem 2.4) that for any $f \in L^{\fhi}(\R)$, there exists $\la>0$ such that:
\[\lim_{\delta \rightarrow 0}\omega(\la f, \delta)_{\fhi}=0,\]
as usually happens in these cases.
Now, in the following, we recall the kind of sampling operators (introduced in \cite{Coroianu}) for which convergence results in $L_{+}^{\fhi}(\Om)$, where $\fhi$ is a convex $\fhi$-function, have been proved in \cite{Boccali}. To do this, first, we recall the definition and the main properties of the kernel functions associated with the max-product Kantorovich sampling operators. From now on, a function $\chi: \R \rightarrow \R$  will be called a generalized kernel if it satisfies the following properties: \vskip0.2cm $(\chi1)$ $\chi \in L^{1}(\R)$ and it is bounded on $\R$; 
\vskip0.1cm $(\chi2)$ For a suitable positive constant $a_{\chi}>0$, we have:\[\inf_{x \in [-1/2,1/2]}\chi(x):=a_{\chi}>0, \textnormal{ if } \Om=\R,\]
\[\hskip0.4cm \inf_{x \in [-3/2,3/2]} \chi(x):=a_{\chi}>0, \ \textnormal{if} \ \Om=[a,b];\] 
\vskip0.1cm $(\chi3)$ for some $\beta >0$, we have: 
\[m_{\beta}(\chi):=\sup_{x \in \R} \bigvee_{k \in \Z} \assolutol \chi(x-k) \assolutor \cdot \assolutol x-k \assolutor^{\beta} < +\infty,\] \hskip0.7cm i.e., the generalized absolute moment of order $\beta$ of $\chi$ is finite. As it is known from the literature (see, e.g., \cite{Coroianu3,Yang,Bede,Anastassio,Gungor,Kadak}), the symbol $\bigvee$ is defined as follows:
\[\bigvee_{k \in \mathcal{J}}A_{k}:= \sup\{A_{k} \in \R, k \in \mathcal{J}\},\] for any set of indexes $\mathcal{J} \subseteq \Z$. It immediately follows that, if $\mathcal{J}$ is a finite set, then $\bigvee$ reduces to a maximum. \newline Now, we recall some useful results concerning the kernel $\chi$ with respect to the symbol $\bigvee$. 
\begin{lemma} [\cite{CCGV}, \textnormal{Lemma 2.1}]
\label{lemma1.1}
Let $\chi$ be a fixed generalized kernel. Then it turns out that:
\[m_{v}(\chi)<+\infty, \quad \textit{for every} \ \ 0 \le v \le \beta,\]
where $\beta >0$ is the constant arising from condition $(\chi3)$. In particular, we have $\mom \le \norma \chi \norma_{\infty}$. 
\end{lemma}
\begin{lemma} [\cite{Coroianu}, \textnormal{Lemma 2.3}]
\label{lemma1.2}
Let $\chi$ be a fixed generalized kernel satisfying $(\chi2)$ with $a_{\chi}>0$. Then, if $\Om=\R$, the following inequality holds:
\begin{equation}
\label{eq2}
\bigvee_{k \in \Z} \chi(nx-k) \ge a_{\chi} >0, \quad x \in \R, \ n \in \N.
\end{equation}
Moreover, if $\Om=[a,b]$, then for all $x \in [a,b]$ there holds:
\begin{equation}
\label{eq3}
\bigvee_{k \in \mathcal{J}_{n}} \chi(nx-k) \ge a_{\chi}>0,
\end{equation}
for every $n \in \N$ sufficiently large, where, here, $\mathcal{J}_{n}:=\{k \in \Z: \lceil na \rceil \le k \le \lfloor nb \rfloor -1\}$. The symbols $\lceil \cdot \rceil$ and $\lfloor \cdot \rfloor$ denote, respectively, the ``ceiling" and the ``integral part" of a given number. 
\end{lemma} 
Now, we are able to recall the definition of max-product Kantorovich sampling operators. For any generalized kernel $\chi$, the corresponding family of max-product Kantorovich sampling operators is defined by:
\[\op(f)(x):=\dfrac{\displaystyle \bigvee_{k \in \mathcal{J}_{n}}\chi(nx-k) \left[n \integralek f(t) \ dt\right]}{\displaystyle \bigvee_{k \in \mathcal{J}_{n}} \chi(nx-k)}, \quad x \in \Om,\]
where $f: \Om \rightarrow \R$ is a locally integrable function. Here, $\mathcal{J}_{n}=\Z$, if $\Om=\R$ or $\mathcal{J}_{n}=\{k \in \Z: \lceil na \rceil \le k \le \lfloor nb \rfloor -1\}$, in the compact case $\Om=[a,b]$. 
\begin{remark}
$(a)$ Keeping assumptions $(\chi1)$ and $(\chi3)$, one can directly assume that for $\chi$ both inequalities (\ref{eq2}) and (\ref{eq3}) hold, instead of assuming  condition $(\chi2)$. It means that we can also consider kernels that are not necessarily bounded from below by a strictly positive parameter in a neighborhood of the origin.
\newline $(b)$ We note that, if $f \in L^{\infty}(\Om)$, by Lemma \ref{lemma1.1} and Lemma \ref{lemma1.2}, we have:
\[\assolutol \op(f)(x) \assolutor \le \dfrac{\mom}{a_{\chi}}\norma f \norma_{\infty}< +\infty,\]
for every $x \in \Om$ and $n \in \N$ sufficiently large, i.e., $\op$ are well-defined and $\op: L^{\infty}(\Om) \rightarrow L^{\infty}(\Om)$ are bounded operators. 
\end{remark}
From the definition of the symbol $\bigvee$, the following useful properties for $\op$ can be easily proved.
\begin{lemma}[\cite{Coroianu}, \textnormal{Lemma 2.6}]
\label{lemma1.3}
Let $\chi$ be a fixed generalized kernel. Then for all $n \in \N$, there holds:
\vskip0.1cm $(i)$ if $f \le g$, then $\op(f) \le \op(g)$;
\vskip0.1cm $(ii)$ $\op(f+g)\le \op(f)+\op(g)$, i.e., $\op$ is sub-additive;
\vskip0.1cm $(iii)$ $\assolutol \op(f)-\op(g) \assolutor \le \op(\assolutol f-g \assolutor)$;
\vskip0.1cm $(iv)$ $\op(\la f)= \la \op(f)$ for each $\la \ge 0$, i.e., $\op$ is positive homogeneous,
\vskip0.2cm \hskip-0.5cm  for every non-negative and bounded functions $f$, $g: \Om \rightarrow \R_{0}^{+}$. 
\end{lemma} 
Now, we recall the following useful inequality concerning any convex $\fhi$-function $\fhi$ in the max-product setting. For a proof, see \cite{Boccali}, Proposition 3.1. 
\begin{lemma}
\label{lemma1.4}
Let $\fhi$ be a convex $\fhi$-function. Then the following inequality:
\begin{equation}
\label{disuguaglianzafhi}
\fhi\left(\bigvee_{k \in \mathcal{J}} A_{k}\right) \le \bigvee_{k \in \mathcal{J}} \fhi(2A_{k}),
\end{equation}
holds for any set of indexes $\mathcal{J} \subseteq \Z$, and $A_{k} \ge 0$, $k \in \mathcal{J}$. 
\end{lemma}
Another preliminary result that will be useful in the next section is the following modular inequality involving $\op$, which has been proved in \cite{Boccali}.
\begin{teorema}[\cite{Boccali}, \textnormal{Theorem 4.1}]
\label{teorema1.1}
Let $\fhi$ be a convex $\fhi$-function and let $\chi$ be a fixed generalized kernel. Then for every $f$, $g \in L_{+}^{\fhi}(\Om)$ and $\la >0$, it turns out that:
\begin{equation}
\label{modularinequality}
\I\bigl[\la\bigl(\op(f)-\op(g)\bigr)\bigr] \le \dfrac{\norma \chi \norma_{1}}{\mom} \I\left[\dfrac{\mom}{a_{\chi}}2\la(f-g)\right],
\end{equation}
for $n \in \N$ sufficiently large. 
\end{teorema}
\begin{remark}
From the above inequality, it follows that $\op$ maps $L_{+}^{\fhi}(\Om)$ into $L^{\fhi}(\Om)$, where $\fhi$ is any convex $\fhi$-function. 
\end{remark}
\section{Main results}
\label{risultatiprincipali}
In this section, we study the order of approximation for the max-product Kantorovich sampling operators, when reconstruction processes of functions belonging to $L^{\fhi}_{+}(\Om)$, with $\fhi$ convex, are considered. In case of $\Om=\R$, in \cite{Coroianu}, Jackson-type estimates with respect to both the sup-norm and $L^{p}$-norm for the family $(\op)_{n \in \N}$ have been established  employing the suitable moduli of smoothness of the (non-negative) function to be approximated. \newline In this paper, we prove a quantitative estimate for the involved operators in the more general setting of Orlicz spaces via the notion of modulus of smoothness recalled in the previous section. In order to reach this purpose, we need to require that the kernel $\chi$ from the definition of $\op$ satisfies the following additional assumption:
\vskip0.1cm $(\chi4)$ for any fixed $0<\alpha<1$, we have:
\[ n \int_{\lvert y \rvert > 1/n^{\alpha}} \assolutol \chi(ny) \assolutor \ dy \le M n^{-\gamma}, \quad \textnormal{for} \ n \in \N \ \textnormal{sufficiently large},\]
and suitable positive absolute constants $M$, $\gamma$ depending on $\alpha$ and $\chi$, and independent on $n \in \N$. 
\begin{teorema}
\label{teorema2.1}
Let $\fhi$ be a convex $\fhi$-function and let $\chi$ be a generalized kernel satisfying condition $(\chi4)$. Moreover, let us denote by $\tau$ the characteristic function of the set $[0,1]$. Then for every $f \in L_{+}^{\fhi}(\R)$ and $\la >0$, it turns out that:
\begin{equation}
\label{assertion}
\begin{split}
\I\bigl[\la\bigl(\op(f)-f\bigr)\bigr]& \le \dfrac{\norma \chi \norma_{1}M_{0}(\tau)}{2 \norma \chi \norma_{\infty}}\omega\left(\dfrac{\norma \chi \norma_{\infty}}{a_{\chi}}4 \la f, \dfrac{1}{n^{\alpha}}\right)_{\fhi} \\ 
& +\frac{M_{0}(\tau)M}{2\norma \chi \norma_{\infty}}I^{\fhi}\left[\frac{\norma \chi \norma_{\infty}}{a_{\chi}}8 \la f\right] \frac{1}{n^{\gamma}} + \dfrac{1}{2}\omega\left(\dfrac{\norma \chi \norma_{\infty}}{a_{\chi}}4 \la f, \dfrac{1}{n}\right)_{\fhi},
\end{split}
\end{equation}
for every $n \in \N$ sufficiently large, where $M$, $\gamma>0$ are the constants of condition $(\chi4)$, and $M_{0}(\tau):=\sup_{x \in \R}\sum_{k \in \Z}\tau(x-k)<+\infty$, since $\tau$ is bounded and with compact support. In particular, for $\la >0$ sufficiently small, the above inequality implies the modular convergence of the max-product Kantorovich sampling operators $\op$ to $f$.
\begin{proof}
Let $\la >0$ be fixed. In order to estimate the quantity $\I[\la(\op(f)-f)]$, we first consider the following family of auxiliary operators defined as:
\[P_{n}^{\chi}(f)(x):=\dfrac{\displaystyle \Vz \chi(nx-k) \left[n \integralek f(t+x-k/n) \ dt\right]}{\displaystyle \Vz \chi(nx-k)}, \quad x \in \R, \ n \in \N.\]
It is easy to prove that the following inequality:
\begin{equation}
\label{eq5}
\assolutol \op(f)(x)-P_{n}^{\chi}(f)(x) \assolutor \le \dfrac{\displaystyle \Vz \chi(nx-k) \left[n \integralek \assolutol f(t)-f(t+x-k/n)\assolutor \ dt\right]}{\displaystyle \Vz \chi(nx-k)},
\end{equation}
holds for all $x \in \R$ and $n \in \N$. Now, we set $\textbf{1}: \R \rightarrow \R$, $\textbf{1}(x):=1$, $x \in \R$ (the identity function) and $f_{x}: \R \rightarrow \R_{0}^{+}$, $f_{x}(t):=f(x)$, $ t \in \R$, for every $x \in \R$. By properties $(iii)$ and $(iv)$ of Lemma \ref{lemma1.3}, which obviously are also satisfied by the auxiliary operators $\pop$ just introduced, together with the fact that $P_{n}^{\chi}(\textbf{1})=\textbf{1}$, we have:
\begin{equation*}
\assolutol \pop (f)(x) -f(x) \assolutor=\assolutol \pop (f)(x)-f(x)\pop(\textbf{1})(x) \assolutor
\end{equation*}
\begin{equation*}
\hskip2.3cm = \assolutol \pop(f)(x)-\pop(f_{x})(x) \assolutor
\end{equation*}
\begin{equation}
\label{eq6}
\hskip1.1cm \le \pop(\assolutol f - f_{x}\assolutor)(x) \\
\end{equation}
\begin{equation*}
\hskip6.6cm = \dfrac{\displaystyle \Vz \chi(nx-k) \left[n \integralek \assolutol f(t+x-k/n)-f(x) \assolutor \ dt \right]}{\displaystyle \Vz \chi(nx-k)},
\end{equation*}
for any $n \in \N$ and $x \in \R$. Since $\fhi$ is non-decreasing and convex, we can write what follows\textcolor{red}{,} for all $n \in \N$:
\begin{equation*}
\begin{split}
\I&\bigl[\la\bigl(\op(f)-f\bigr)\bigr] \\
& \le \frac{1}{2} \left\{\integraler \fhi\bigl(2\la\assolutol\op(f)(x)-\pop(f)(x)\assolutor\bigr) \ dx + \integraler \fhi\bigl(2\la\assolutol\pop(f)(x)-f(x)\assolutor\bigr) \ dx\right\} \\
&=:\frac{1}{2}\bigl\{I_{1}+I_{2}\bigr\}.
\end{split}
\end{equation*}
Now, by using  (\ref{eq5}), (\ref{eq2}), the inequality recalled in Lemma \ref{lemma1.4}, and the convexity of $\fhi$ with the fact that $\assolutol \chi(nx-k) \assolutor \le \norma \chi \norma_{\infty}$, for every $k \in \Z$, we obtain: 
\begin{equation*}
\begin{split}
I_{1}&\le \integraler \fhi\left(\frac{2\la}{a_{\chi}}\Vz \assolutol\chi(nx-k) \assolutor \left[n \integralek \assolutol f(t)-f(t+x-k/n)\assolutor \ dt\right]\right) \ dx \\
&\le \integraler \Vz \fhi\left(\frac{4\la}{a_{\chi}}\assolutol\chi(nx-k)\assolutor \left[n \integralek \assolutol f(t)-f(t+x-k/n)\assolutor \ dt\right]\right) \ dx \\
& \le \integraler \Vz \dfrac{\assolutol\chi(nx-k)\assolutor}{\norma \chi \norma_{\infty}} \fhi\left(\frac{\norma \chi\norma_{\infty}}{a_{\chi}}4\la\left[n \integralek  \assolutol f(t)-f(t+x-k/n)\assolutor \ dt\right]\right) \ dx.
\end{split}
\end{equation*}
Moreover, we claim that the following inequality can be stated:
\begin{equation}
\label{somma}
\Vz A_{k} \le \sum_{k \in \Z} A_{k},
\end{equation}
for any $A_{k} \ge 0$, $k \in \Z$. Indeed, from the definition of the symbol $\bigvee$, it follows that for every $\varepsilon >0$ there exists $\overline{k} \in \Z$ such that:
\[\Vz A_{k}<A_{\overline{k}}+\varepsilon<\sum_{k \in \Z} A_{k}+\varepsilon.\]
Hence, the assertion follows by the arbitrariness of $\varepsilon >0$. Now, by using the Jensen inequality (see, e.g., \cite{Jensen}), the inequality (\ref{somma}), the Fubini-Tonelli theorem, and the change of variable $y=x-k/n$, we have:
\begin{equation*}
\begin{split}
I_{1}& \le \integraler \Vz \dfrac{\assolutol \chi(nx-k)\assolutor}{\norma \chi \norma_{\infty}} \left[n \integralek \fhi\left(\dfrac{\norma \chi \norma_{\infty}}{a_{\chi}}4\la\assolutol f(t)-f(t+x-k/n)\assolutor\right) \ dt\right] \ dx \\
&\le \integraler \sum_{k \in \Z} \dfrac{\assolutol \chi(nx-k)\assolutor}{\norma \chi \norma_{\infty}} \left[n \integralek \fhi\left(\dfrac{\norma \chi \norma_{\infty}}{a_{\chi}}4\la\assolutol f(t)-f(t+x-k/n) \assolutor\right) \ dt\right] \ dx \\
&= \sum_{k \in \Z} \integraler \dfrac{\assolutol \chi(nx-k)\assolutor}{\norma \chi \norma_{\infty}} \left[n \integralek \fhi\left(\dfrac{\norma \chi \norma_{\infty}}{a_{\chi}}4\la\assolutol f(t)-f(t+x-k/n) \assolutor\right) \ dt\right] \ dx \\
& = \integraler  \dfrac{\assolutol \chi(ny)\assolutor}{\norma \chi \norma_{\infty}} \sum_{k \in \Z} \left[n \integralek \fhi\left(\dfrac{\norma \chi \norma_{\infty}}{a_{\chi}}4\la \assolutol f(t)-f(t+y) \assolutor\right) \ dt\right] \ dy \\
&= \integraler  \dfrac{\assolutol \chi(ny)\assolutor}{\norma \chi \norma_{\infty}} \sum_{k \in \Z} \left[n \integraler \fhi\left(\dfrac{\norma \chi \norma_{\infty}}{a_{\chi}}4\la \assolutol f(t)-f(t+y) \assolutor\right) \tau(nt-k) \ dt\right] \ dy \\
&= \integraler \dfrac{\assolutol \chi(ny) \assolutor}{\norma \chi \norma_{\infty}} \left[n \integraler\fhi\left(\dfrac{\norma \chi \norma_{\infty} }{a_{\chi}}4\la \assolutol f(t)-f(t+y)\assolutor\right) \sum_{k \in \Z} \tau(nt-k) \ dt\right] \ dy \\
&\le M_{0}(\tau) \norma \chi \norma_{\infty}^{-1} \integraler n \assolutol \chi(ny) \assolutor \left[\integraler \fhi\left(\dfrac{\norma \chi \norma_{\infty}}{a_{\chi}}4\la \assolutol f(t)-f(t+y)\assolutor\right) \ dt \right] \ dy \\
&=M_{0}(\tau)\norma \chi \norma_{\infty}^{-1} \integraler n \assolutol \chi(ny) \assolutor \I\left[\dfrac{\norma \chi \norma_{\infty}}{a_{\chi}}4\la \bigl(f(\cdot)-f(\cdot+y)\bigr)\right] \ dy=: J,
\end{split}
\end{equation*}
where $M_{0}(\tau)$ denotes the usual discrete absolute moment of order $0$ of $\tau$, where $\tau$ is the characteristic function of the set $[0,1]$. Now, let $0<\alpha<1$ be fixed. Thus, we can split the above integral $J$ as follows:
\begin{equation*}
\begin{split}
J= M_{0}(\tau)\norma \chi \norma_{\infty}^{-1} \left\{\int_{\lvert y \rvert\le 1/n^{\alpha}}+\int_{\lvert y \rvert > 1/n^{\alpha}}\right\} n \assolutol \chi(ny) \assolutor \I\left[\dfrac{\norma \chi \norma_{\infty}}{a_{\chi}}4\la \bigl( f(\cdot)-f(\cdot+y)\bigr)\right] \ dy=: J_{1}+J_{2}.
\end{split}
\end{equation*}
Estimating $J_{1}$, we obtain:
\begin{equation*}
\begin{split}
J_{1}&\le M_{0}(\tau) \norma \chi \norma_{\infty}^{-1} \int_{\lvert y \rvert\le 1/n^{\alpha}} n \assolutol \chi(ny) \assolutor \omega\left(\dfrac{\norma \chi \norma_{\infty}}{a_{\chi}}4\la f, \lvert y \rvert\right)_{\fhi} \ dy  \\
&\le M_{0}(\tau)\norma \chi \norma_{\infty}^{-1} \omega\left(\dfrac{\norma\chi \norma_{\infty}}{a_{\chi}}4\la f, \dfrac{1}{n^{\alpha}}\right)_{\fhi} \int_{\lvert y \rvert\le 1/n^{\alpha}} n \assolutol \chi(ny) \assolutor \ dy \\
&\le \dfrac{\norma \chi \norma_{1}M_{0}(\tau)}{\norma \chi \norma_{\infty}} \omega\left(\dfrac{\norma \chi \norma_{\infty}}{a_{\chi}}4\la f, \dfrac{1}{n^{\alpha}}\right)_{\fhi},
\end{split}
\end{equation*}
for every $n \in \N$. Concerning $J_{2}$, since $\fhi$ is convex, one has:
\begin{equation*}
\begin{split}
J_{2}& \le  M_{0}(\tau)\norma \chi \norma_{\infty}^{-1} \int_{\lvert y \rvert > 1/n^{\alpha}} n \assolutol \chi(ny) \assolutor \frac{1}{2} \left\{\I\left[\frac{\norma \chi \norma_{\infty}}{a_{\chi}}8\la f\right]+ \I\left[\dfrac{\norma \chi \norma_{\infty}}{a_{\chi}}8\la f(\cdot+y)\right]\right\} \ dy.
\end{split}
\end{equation*}
Thus, observing that:
\[\I\left[\dfrac{\norma \chi \norma_{\infty}}{a_{\chi}}8\la f\right]=\I\left[\dfrac{\norma \chi \norma_{\infty}}{a_{\chi}}8\la f(\cdot+y)\right],\]
for every $y \in \R$, and since condition $(\chi4)$ holds, we finally get:
\begin{equation*}
\begin{split}
J_{2}&\le  M_{0}(\tau)\norma \chi \norma_{\infty}^{-1} \I\left[\dfrac{\norma \chi \norma_{\infty}}{a_{\chi}}8\la f\right] n \int_{\lvert y \rvert > 1/n^{\alpha}} \assolutol \chi(ny) \assolutor \ dy \\
& \le \dfrac{M_{0}(\tau)}{\norma \chi \norma_{\infty}} \I\left[\dfrac{\norma \chi \norma_{\infty}}{a_{\chi}}8\la f\right] Mn^{-\gamma},
\end{split}
\end{equation*}
for $n \in \N$ sufficiently large. Now, let us estimate $I_{2}$.  By using (\ref{eq6}), and proceeding as in the previous estimates, we can write what follows:
\begin{equation*}
\begin{split}
I_{2} &\le \integraler \fhi\left(\dfrac{2\la}{a_{\chi}}\Vz \assolutol \chi(nx-k)\assolutor \left[n \integralek \assolutol f(t+x-k/n)-f(x)\assolutor \ dt\right]\right) \ dx \\
& \le \integraler \Vz \fhi\left(\dfrac{4\la}{a_{\chi}}\assolutol \chi(nx-k) \assolutor \left[n \integralek \assolutol f(t+x-k/n)-f(x) \assolutor \ dt\right]\right) \ dx \\
&\le \integraler \Vz \dfrac{\assolutol \chi(nx-k)\assolutor}{\norma \chi \norma_{\infty}} \fhi\left(\dfrac{\norma \chi \norma_{\infty}}{a_{\chi}}4\la\left[n \integralek \assolutol f(t+x-k/n)-f(x)\assolutor \ dt\right]\right) \ dx,
\end{split}
\end{equation*}
for any $n \in \N$. Exploiting the Jensen inequality, as above, the change of variable $y=t-k/n$, and the Fubini-Tonelli theorem, we finally obtain:
\begin{equation*}
\begin{split}
I_{2}&\le \norma \chi \norma_{\infty}^{-1} \integraler \Vz \assolutol \chi(nx-k) \assolutor \left[n \integralek \fhi\left(\dfrac{\norma \chi \norma_{\infty}}{a_{\chi}}4\la \assolutol f(t+x-k/n)-f(x) \assolutor\right) \ dt \right] \ dx \\
& = \norma \chi \norma_{\infty}^{-1} \integraler \Vz \assolutol \chi(nx-k) \assolutor \left[n \int_{0}^{1/n}\fhi\left(\dfrac{\norma \chi \norma_{\infty}}{a_{\chi}}4\la \assolutol f(x+y)-f(x) \assolutor\right) \ dy\right] \ dx \\
&= \norma \chi \norma_{\infty}^{-1} \integraler \left[ n \int_{0}^{1/n} \fhi\left(\dfrac{\norma \chi \norma_{\infty}}{a_{\chi}}4 \la \assolutol f(x+y)-f(x)\assolutor\right) \ dy \right] \Vz \assolutol \chi(nx-k) \assolutor \ dx \\
& \le \integraler n \int_{0}^{1/n} \fhi\left(\dfrac{\norma \chi\norma_{\infty}}{a_{\chi}}4\la \assolutol f(x+y)-f(x) \assolutor\right) \ dy \ dx \\ 
&=n \int_{0}^{1/n} \I\left[\dfrac{\norma \chi \norma_{\infty}}{a_{\chi}}4 \la \bigl(f(\cdot+y)-f(\cdot)\bigr)\right] \ dy \\
& \le n \int_{0}^{1/n} \omega\left(\dfrac{\norma \chi \norma_{\infty}}{a_{\chi}}4\la f, y\right)_{\fhi} \ dy \le \omega\left(\dfrac{\norma \chi \norma_{\infty}}{a_{\chi}}4\la f, \dfrac{1}{n}\right)_{\fhi},
\end{split}
\end{equation*}
for all $n \in \N$. This completes the proof. 
\end{proof}  
\end{teorema} 
\begin{remark}
$(a)$ Note that, it is easy to show that for any kernel $\chi$ such that $\chi(x)=\mathcal{O}(\lvert x\rvert^{-\theta})$, as $\lvert x \rvert \rightarrow +\infty$, for $\theta >1$, we have that assumption $(\chi4)$ is satisfied for some constant $M>0$ and $\gamma=(1-\alpha)(\theta-1)>0$, for every fixed $0<\alpha<1$ (see, e.g., \cite{Costarelli}). 
\newline $(b)$ For some examples of generalized kernels, both with and without compact support, which satisfy condition $(\chi4)$, one can see \cite{Coroianu}. 
\end{remark}
\begin{remark}
The estimate established in Theorem \ref{teorema2.1} can be extended to the case of functions of arbitrary sign that are bounded from below. Indeed, let $f: \R \rightarrow \R$ be fixed such that $f(x)\ge c$, for every $x \in \R$, with $c \in \R$. Then, it is easy to see that for the family of max-product operators defined by $\overline{K}_{n}^{\chi}(f):=\op(f(\cdot)-c)+c$, assuming $f$ in suitable spaces, it turns out that:
\begin{equation*}
\begin{split}
\I\bigl[\la\bigl(\overline{K}_{n}^{\chi}(f)-f\bigr)\bigr] & \le \dfrac{\norma \chi \norma_{1}M_{0}(\tau)}{2 \norma \chi \norma_{\infty}}\omega\left(\dfrac{\norma \chi \norma_{\infty}}{a_{\chi}}4 \la f, \dfrac{1}{n^{\alpha}}\right)_{\fhi} \\ 
& +\frac{M_{0}(\tau)M}{2\norma \chi \norma_{\infty}}I^{\fhi}\left[\frac{\norma \chi \norma_{\infty}}{a_{\chi}}8 \la (f(\cdot)-c)\right] \frac{1}{n^{\gamma}} + \dfrac{1}{2}\omega\left(\dfrac{\norma \chi \norma_{\infty}}{a_{\chi}}4 \la f, \dfrac{1}{n}\right)_{\fhi},  
\end{split}
\end{equation*} 
for every $n \in \N$ sufficiently large.  
\end{remark}
From the quantitative estimate (\ref{assertion}) of Theorem \ref{teorema2.1}, it is possible to determine the qualitative order of approximation, assuming $f:\R \rightarrow \R$ in suitable Lipschitz classes. More precisely, we recall that, for any fixed $0 < v \le 1$, the Lipschitz class in Orlicz spaces $L^{\fhi}(\R)$, denoted by $Lip_{\fhi}(v)$, can be defined as the space of all functions $f \in L^{\fhi}(\R)$ such that there exists $\la >0$ for which:
\[\I\bigl[\la\bigl(f(\cdot+h)-f(\cdot)\bigr)\bigr]=\integraler \fhi\bigl(\la\assolutol f(x+h)-f(x) \assolutor\bigr) \ dx=\mathcal{O}(\lvert h \rvert^{v}),\]
as $\lvert h \rvert \rightarrow 0$.  As a direct consequence of Theorem \ref{teorema2.1}, we can state the following.
\begin{cor}
Under the assumptions of Theorem \ref{teorema2.1} with $0 < \alpha <1$, for any non-negative $f \in Lip_{\fhi}(v)$, with $0 < v \le 1$, there exist $K>0$ and $\la >0$ such that: 
\[\I\bigl[\la\bigl(\op(f)-f\bigr)\bigr] \le Kn^{-\rho},\]
for $n \in \N$ sufficiently large, where $\rho:=\min\{\alpha v,\gamma\}$ and $\gamma$ is the constant from condition $(\chi4)$. 
\end{cor} 
In the remaining part of the paper, we consider the approximation of functions defined on compact intervals $[a,b]$. In \cite{Coroianu}, quantitative properties  for the operators $\op$ with respect to the norm $\norma \cdot \norma_{p}$, $1\le p < +\infty$, have been studied using the well-known definition of the K-functionals introduced by Peetre (see, e.g., \cite{Peetre}). In order to investigate the order of approximation for the max-product Kantorovich neural network (NN) operators in $L_{+}^{\fhi}([a,b])$, a suitable notion of K-functional has been introduced in \cite{CostaSambu}. Here, we exploit such definition in order to reach the same purpose for $\op$. Note that, in the case of functions defined on $[a,b]$, in order to obtain quantitative estimates for the order of approximation, it is preferable to work with K-functionals in place of moduli of smoothness, due to technical reasons. \newline Let us denote by $C_{+}^{1}([a,b])$ the space of non-negative differentiable functions with continuous derivative. For any fixed $f \in L_{+}^{\fhi}([a,b])$, we define the Orlicz-type Peetre K-functional as follows:
\begin{equation}
\mathcal{K}(f,\la,\delta)_{\fhi}:= \inf_{g \in C_{+}^{1}([a,b])} \{\I\bigl[\la\bigl(f-g\bigr)\bigr]+\delta\fhi(\norma g'\norma_{\infty})\}, \quad \delta>0,
\end{equation}
for some $\la >0$. It is easy to see that $C_{+}^{1}([a,b]) \subset L_{+}^{\fhi}([a,b])$. Moreover, we can prove that for any fixed $f \in L^{\fhi}_{+}([a,b])$, there exists $\la >0$ such that:
\[\mathcal{K}(f,\la,\delta)_{\fhi} \le \I[2\la f]+\I[2\la g]+\delta\fhi(\norma g'\norma_{\infty})<+\infty,\]
for every $g \in C_{1}^{+}([a,b])$. This means that the definition of $\mathcal{K}(f, \la, \delta)_{\fhi}$ is well-posed. The above quantity provides a useful tool for investigating the approximation properties and the smoothness of $f$. Indeed, if the inequality $\mathcal{K}(f, \la, \delta)_{\fhi} < \varepsilon$ holds for some $\la >0$ and $\delta >0$, then there exists $g \in C_{1}^{+}([a,b])$, which approximates $f$ with an error $\I[\la(f-g)] < \varepsilon$, and such that $\fhi(\norma g' \norma_{\infty})<\delta^{-1}\varepsilon$, i.e, $\norma g' \norma_{\infty}$ is not too large. For more details concerning general definitions of K-functionals, see \cite{Vore}, Chapter 6. \newline Now, we are able to give the proof of the main theorem. 
\begin{teorema}
\label{teorema2.2}
Let $\fhi$ be a convex $\fhi$-function and let $\chi$ be a generalized kernel which satisfies $(\chi3)$ with $\beta \ge 1$. Then for any $f \in L_{+}^{\fhi}([a,b])$, there exist two suitable positive constants $\la_{0}$, $\la_{1}>0$ such that:
\begin{equation}
\label{keyfunctional}
\I\bigl[\la_{1}\bigl(\op(f)-f\bigr)\bigr] \le A_{1}\mathcal{K}\left(f,\la_{0}, A_{2}\dfrac{1}{n}\right)_{\fhi},
\end{equation}
for every $n \in \N$ sufficiently large, where: 
\[A_{1}:= \dfrac{\norma \chi \norma_{1}}{\mom}+1, \quad \textnormal{and} \quad A_{2}:=\dfrac{\la_{0}\left(\frac{1}{2}\mom+m_{1}(\chi)\right)(b-a)}{a_{\chi}\left(\dfrac{\norma \chi \norma_{1}}{\mom}+1\right)}.\]
Note that $\mom$ and $m_{1}(\chi)$ are both finite in view of Lemma \ref{lemma1.1}. 
\begin{proof}
First of all, in view of the definition of the space $L_{+}^{\fhi}([a,b])$, there exists $\la_{0} >0$ such that $\I[\la_{0}f]< +\infty$. Now, let $g \in C_{+}^{1}([a,b])$ be fixed and let us choose $\la_{1} >0$ such that:
\[\dfrac{\mom}{a_{\chi}}6\la_{1}\le \la_{0}.\] 
Then, by using the properties of the modular functional $\I$, the inequality of Theorem \ref{teorema1.1}, and observing that $\mom \ge a_{\chi}$, we can write what follows:
\begin{equation}
\label{modulare}
\begin{split}
\I\bigl[\la_{1}\bigl(\op(f)-f\bigr)\bigr]& \le \frac{1}{3}\bigl\{\I\bigl[3\la_{1}\bigl(\op(f)-\op(g)\bigr)\bigr]+\I\bigl[3\la_{1}\bigl(\op(g)-g\bigr)\bigr]+\I\bigl[3\la_{1}\bigl(g-f\bigr)\bigr]\bigr\} \\
&\le \dfrac{\norma \chi \norma_{1}}{\mom}\I\left[\dfrac{\mom}{a_{\chi}}6\la_{1}(f-g)\right] + \I\bigl[3\la_{1}\bigl(\op(g)-g\bigr)\bigr]+ \I\bigl[3\la_{1}\bigl(f-g\bigr)\bigr] \\
&\le \left(\dfrac{\norma \chi \norma_{1}}{\mom}+1\right)\I\bigl[\la_{0}\bigl(f-g\bigr)\bigr] + \I\bigl[\la_{0}\bigl(\op(g)-g\bigr)\bigr],
\end{split}
\end{equation}
for $n \in \N$ sufficiently large. Now, for every fixed $x \in [a,b]$, we consider the function $g_{x}:[a,b]\rightarrow\R$, $g_{x}(t)=g(x)$, $t \in [a,b]$. By the properties of $\op$ established in Lemma \ref{lemma1.3}, proceeding as in the proof of Theorem \ref{teorema2.1}, we have:
\begin{equation*}
\begin{split}
\assolutol \op(g)(x)-g(x)\assolutor&=\assolutol \op(g)(x)-g(x)\op(\textbf{1})(x) \assolutor \\
&=\assolutol \op(g)(x)-\op(g_{x})(x) \assolutor \\
& \le \op(\assolutol g-g_{x}\assolutor)(x).
\end{split}
\end{equation*}
In the previous estimate, we used the fact that $\op(\textbf{1})=\textbf{1}$, where $\textbf{1}(x)$ denotes again the identity function. Since $g \in C_{+}^{1}([a,b])$, it is easy to prove that:
\[\assolutol g(t)-g_{x}(t)\assolutor=\assolutol g(t)-g(x)\assolutor\le \norma g'\norma_{\infty} \assolutol t-x \assolutor=:\norma g' \norma_{\infty} \eta_{x}(t),\]
for every $t \in [a,b]$. Thus, by $(i)$ and $(iv)$ of Lemma \ref{lemma1.3}, and using (\ref{eq3}), we obtain for $n \in \N$ sufficiently large: 
\begin{equation*}
\begin{split}
\op(\assolutol g-g_{x}\assolutor)(x) &\le\norma g' \norma_{\infty}\op(\eta_{x})(x) \\
&=\norma g' \norma_{\infty} \dfrac{\displaystyle \Vj \chi(nx-k) \left[n \integralek \assolutol t-x \assolutor \ dt\right] }{\displaystyle \Vj \chi(nx-k)} \\
&\le \dfrac{\norma g' \norma_{\infty}}{a_{\chi}} \Vj \assolutol \chi(nx-k) \assolutor \left[n \integralek \assolutol t-x \assolutor \ dt\right] \\
&\le \dfrac{\norma g' \norma_{\infty}}{a_{\chi}} \Vj \assolutol \chi(nx-k) \assolutor \left[ n \integralek \assolutol t-k/n \assolutor \ dt\right]
\end{split}
\end{equation*}
\begin{equation*}
\begin{split}
\hskip3.6cm &+ \dfrac{\norma g' \norma_{\infty}}{a_{\chi}} \Vj \assolutol\chi(nx-k) \assolutor \left[n \integralek \assolutol x-k/n \assolutor \ dt \right] \\
&= \dfrac{\norma g' \norma_{\infty}}{na_{\chi}} \left\{\frac{1}{2}\Vj \assolutol \chi(nx-k) \assolutor +\Vj \assolutol \chi(nx-k) \assolutor \assolutol nx-k \assolutor\right\} \\
&\le \dfrac{\norma g' \norma_{\infty}}{na_{\chi}}\left(\frac{1}{2}\mom+m_{1}(\chi)\right)<+\infty,
\end{split}
\end{equation*}
in view of Lemma \ref{lemma1.1}, since $\beta \ge 1$. Hence, we proved that:
\[\assolutol \op(g)(x)-g(x) \assolutor \le \dfrac{\norma g' \norma_{\infty}}{na_{\chi}}\left(\frac{1}{2}\mom+m_{1}(\chi)\right),\]
for every $x \in [a,b]$ and $n \in \N$ sufficiently large. Then, since $\fhi$ is convex, we get:
\begin{equation*}
\begin{split}
\I\bigl[\la_{0}\bigl(\op(g)-g\bigr)\bigr] &\le \integraleab \fhi\left(\la_{0}\dfrac{\norma g' \norma_{\infty}}{na_{\chi}}\left(\frac{1}{2}\mom+m_{1}(\chi)\right)\right) \ dx \\
&=\fhi\left(\la_{0}\dfrac{\norma g' \norma_{\infty}}{na_{\chi}}\left(\frac{1}{2}\mom+m_{1}(\chi)\right)\right) (b-a) \\
&\le \dfrac{\la_{0}}{na_{\chi}}\left(\frac{1}{2}\mom+m_{1}(\chi)\right)\fhi\bigl(\norma g'\norma_{\infty}\bigr)(b-a),
\end{split}
\end{equation*} 
for every $n \in \N$ sufficiently large. Hence, using this result in (\ref{modulare}), we finally obtain:
\begin{equation*}
\begin{split}
\I&\bigl[\la_{1}\bigl(\op(f)-f\bigr)\bigr] \\ &\le \left(\dfrac{\norma \chi \norma_{1}}{\mom}+1\right)\I\bigl[\la_{0}\bigl(f-g\bigr)\bigr]+\dfrac{\la_{0}}{na_{\chi}}\left(\frac{1}{2}\mom+m_{1}(\chi)\right)(b-a)\fhi\bigl(\norma g'\norma_{\infty}\bigr).
\end{split}
\end{equation*}
for every $n \in \N$ sufficiently large. Thus, the proof follows by the arbitrariness of $g \in C_{+}^{1}([a,b])$. 
\end{proof}
\end{teorema}    
\begin{remark}
$(a)$ Note that, if we consider functions that are not-necessarily non-negative and bounded from below, the estimate established in Theorem \ref{teorema2.2} can be restated. Indeed, let $f:[a,b] \rightarrow \R$ be fixed. Denoting $c:=\inf_{x \in [a,b]} f(x)$, it is easy to see that for the family of max-product operators defined by $\overline{K}_{n}^{\chi}(f):=\op(f(\cdot)-c)+c$, we have:
\[\I\bigl[\la_{1}\bigl(\overline{K}_{n}^{\chi}(f)-f\bigr)\bigr]=\I\bigl[\la_{1}\bigl(\op(f-c)-(f-c)\bigr)\bigr].\]
Moreover, the right side of (\ref{keyfunctional}) can be obtained, observing that:
\begin{equation*}
\begin{split}
\mathcal{K}\bigl(f-c,\la, \delta\bigr)_{\fhi}&=\inf_{g \in C_{+}^{1}([a,b])} \bigl\{\I\bigl[\la\bigl((f-c)-g\bigr)\bigr]+\delta\fhi\bigl(\norma g' \norma_{\infty}\bigr)\bigr\} \\
&=\inf_{g \in C_{+}^{1}([a,b])} \bigl\{\I\bigl[\la\bigl(f-(g+c)\bigr)\bigr]+\delta\fhi\bigl(\norma(g+c)'\norma_{\infty}\bigr)\bigr\} \\
&=\inf_{h \in C^{1}_{+}([a,b]), \ h \ge c} \bigl\{\I\bigl[\la\bigl(f-h\bigr)\bigr]+\delta\fhi\bigl(\norma h' \norma_{\infty}\bigr)\bigr\},
\end{split}
\end{equation*}
for every $\la$, $\delta >0$. 
\newline $(b)$ In order to find examples of generalized kernels for which Theorem \ref{teorema2.2} can be applied, we recall the following useful sufficient condition concerning the existence of generalized absolute moments (see \cite{CCGV}, Lemma 2.1). If $\chi:\R \rightarrow \R$ is a bounded function such that $\chi(x)=\mathcal{O}(\lvert x \rvert^{-\alpha})$, as $\lvert x \rvert \rightarrow +\infty$, with $\alpha \ge 1$. Then there holds:
\[m_{\beta}(\chi)< +\infty, \ \  \textnormal{ for every} \ 0 \le \beta \le \alpha.\]
\end{remark}
We want to point out that all the results proved in this section hold for (non-negative) functions that belong to the Orlicz space $L^{\fhi}(\Om)$, with $\fhi$ convex, and where $\Om$ can be (again) a bounded interval or the whole real axis. Therefore, the above results are part of a unifying theory that can be applied to a wide range of functional spaces, including the usual $L^{p}$-spaces, the Zygmund (or interpolation) spaces and the exponential spaces. The first ones can be obtained by choosing $\fhi(u)=u^{p}$, $u \ge 0$, $1 \le p < +\infty$, which implies that $\I[f]=\norma f \norma_{p}^{p}$. The convex $\fhi$-functions $\fhi_{\alpha,\beta}(u):=u^{\alpha}\log^{\beta}(u+e)$, $u \ge 0$, with $\alpha \ge 1$ and $\beta >0$, generate the so-called interpolation spaces, which find applications in the theory of partial differential equations. Other relevant examples of Orlicz spaces are provided by the well-known exponential spaces, generated by $\fhi_{\gamma}(u):=e^{u^{\gamma}}-1$, $u \ge 0$, with $\gamma >0$. The latter are widely used to obtain embedding theorems between Sobolev spaces. The modular functionals corresponding to $\fhi_{\alpha,\beta}$ and $\fhi_{\gamma}$ are:
\[I^{\fhi_{\alpha,\beta}}[f]:=\integrale \assolutol f(x) \assolutor^{\alpha} \log^{\beta}\bigl(\assolutol f(x) \assolutor+e\bigr) \ dx, \quad f \in M(\Om),\] 
and 
\[I^{\fhi_{\gamma}}[f]:=\integrale \left(e^{\lvert f(x) \rvert^{\gamma}}-1\right) \ dx, \quad f \in M(\Om),\]
respectively. For more details about the above functional spaces, one can see, e.g., \cite{Edmund,Rao,Bardaro3,Tarsi}.    


\section*{Acknowledgments}

{\small The authors are members of the Gruppo Nazionale per l'Analisi Matematica, la Probabilit\`a e le loro Applicazioni (GNAMPA) of the Istituto Nazionale di Alta Matematica (INdAM), of the Gruppo UMI (Unione Matematica Italiana) T.A.A. (Teoria dell'Approssimazione e Applicazioni)}, and of the network RITA (Research ITalian network on Approximation).

\section*{Funding}

{\small The second and the third authors have been partially supported within the (1) 2022 GNAMPA-INdAM Project "Enhancement e segmentazione di immagini mediante operatori di tipo campionamento e metodi variazionali per lo studio di applicazioni biomediche'', (2) "Metodiche di Imaging non invasivo mediante angiografia OCT sequenziale per lo studio delle Retinopatie degenerative dell'Anziano (M.I.R.A.)", funded by the Fondazione Cassa di Risparmio di Perugia (FCRP), 2019 and (3) "National Innovation Ecosystem grant ECS00000041 - VITALITY", funded by the European Union - NextGenerationEU under the Italian Ministry of University and Research (MUR). Moreover, the second author has been partially supported within the 2023 GNAMPA-INdAM Project "Approssimazione costruttiva e astratta mediante operatori di tipo sampling e loro applicazioni".
	
}

\section*{Conflict of interest/Competing interests}

{\small The authors declare that they have no conflict of interest and competing interest.}

\section*{Availability of data and material and Code availability}

{ \small Not applicable.}



\begin{thebibliography}{99}
\bibitem{Anastassio}
G. A. Anastassiou, \textit{Nonlinearity: Ordinary and Fractional Approximations by Sublinear and Max-Product Operators}, (Springer, Heidelberg, New York, 2018).

\bibitem{Bardaro2}
C. Bardaro, P. L. Butzer, R. L. Stens and G. Vinti, \textit{Kantorovich-type generalized sampling series in the setting of Orlicz spaces}, Sampl. Theory Signal Image Process., \textbf{6}(1) (2007), 29-52.   	

\bibitem{Bardaro3}
C. Bardaro, J. Musielak and G. Vinti, \textit{Nonlinear Integral Operators and Applications}, De Gruyter Series in Nonlinear Analysis and Applications (New York, Berlin, 2003). 

\bibitem{Bede2}
B. Bede, L. Coroianu and S. G. Gal, \textit{Approximation and shape preserving properties of the Bernstein operator of max-product kind}, Int. J. Math. Math. Sci. (2009). 

\bibitem{Bede3}
B. Bede, L. Coroianu and S. G. Gal, \textit{Approximation and shape preserving properties of the nonlinear Favard-Sz\'asz-Mirakjan operator of max-product kind}, Filomat, \textbf{24}(3) (2010), 55-72. 

\bibitem{Bede4}
B. Bede, L. Coroianu and S. G. Gal, \textit{Approximation of fuzzy numbers by max-product Bernstein operators}, Fuzzy Sets and Systems, \textbf{257} (2014) 41-66.

\bibitem{Bede}
B. Bede, L. Coroianu and S. G. Gal, \textit{Approximation by Max-Product Type Operators}, (Springer, New York, 2016).

\bibitem{Benedetto}
J. J. Benedetto and P. J. S. G. Ferreira, \textit{Modern Sampling Theory: Mathematics and Applications}, (Birkh\"auser, Boston-Basel-Berlin, 2001). 

\bibitem{Boccali}
L. Boccali, D. Costarelli and G. Vinti, \textit{Convergence results in Orlicz spaces for sequences of max-product Kantorovich sampling operators}, ArXiv:https://arxiv.org/abs/2305.18783, submitted (2023).


\bibitem{Butzer}
P. L. Butzer, \textit{A survey of the Whittaker-Shannon sampling theorem and some of its extensions}, J. Math. Res. Exposition, \textbf{3} (1983) 185-212. 

\bibitem{Butzer2}
P. L. Butzer, S. Ries and R. L. Stens, \textit{Approximation of continuous and discontinuous functions by generalized sampling series}, J. Approx. Theory, \textbf{50} (1987) 25-39.

\bibitem{Butzer3}
P. L. Butzer and R. L. Stens, \textit{Reconstruction of signals in $L^{p}(\R)$-space by generalized sampling series based on linear combinations of B-splines}, Integral Transforms Spec. Funct., \textbf{19}(1) (2008) 35-58.   

\bibitem{Cornejo}
M. E. Cornejo, D. Lobo and J. Medina, \textit{On the solvability of bipolar max-product fuzzy relation equations with the standard negation}, Fuzzy Sets and Systems, \textbf{410} (2021) 1-18.

\bibitem{CCGV}
L. Coroianu, D. Costarelli, S. G. Gal and G. Vinti, \textit{The max-product generalized sampling operators: Convergence and quantitative estimates}, Appl. Math. Comput., \textbf{355} (2019) 173-183.  

\bibitem{Coroianu}
L. Coroianu, D. Costarelli, S. G. Gal and G. Vinti, \textit{Approximation by max-product sampling Kantorovich operators with generalized kernels}, Anal. Appl., \textbf{19} (2021) 219-244.

\bibitem{Coroianu3}
L. Coroianu and S. G. Gal, \textit{Approximation by nonlinear generalized sampling operators of max-product kind}, Sampl. Theory Signal Image Process., \textbf{9}(1-3) (2010) 59-75. 

\bibitem{CG2}
L. Coroianu and S. G. Gal, \textit{Localization results for the Bernstein max-product operator}, Appl. Math. Comput., \textbf{231} (2014) 73-78. 

\bibitem{CG}
L. Coroianu and S. G. Gal, \textit{Localization results for the non-truncated max-product sampling operators based on Fejér and sinc-type kernels}, Demonstr. Math., \textbf{49}(1), 38-49 (2016).

\bibitem{Coroianu2}
L. Coroianu and S. G. Gal, \textit{Approximation by max-product operators of Kantorovich type}, Stud. Univ. Babes-Bolyai Math., \textbf{64}(2) (2019) 207-223. 

\bibitem{CostaSambu}
D. Costarelli and A. R. Sambucini, \textit{Approximation results in Orlicz spaces for sequences of Kantorovich max-product neural networks operators}, Results Math., \textbf{73}(1) (2018) 1-15.

\bibitem{CSV}
D. Costarelli, A. R. Sambucini and G. Vinti, \textit{Convergence in Orlicz spaces by means of the multivariate max-product neural network operators of the Kantorovich type and applications}, Neural Comput. Appl., \textbf{31} (2019) 5069-5078.

\bibitem{Jensen}
D. Costarelli and R. Spigler, \textit{How sharp is the Jensen inequality?}, J. Inequal. Appl., \textbf{2015}(69) (2015) 1-10.  

\bibitem{CV}
D. Costarelli and G. Vinti, \textit{Saturation classes for max-product neural network operators activated by sigmoidal functions}, Results Math., \textbf{72}(3) (2017) 1555-1569.

\bibitem{Costarelli}
D. Costarelli and G. Vinti, \textit{A quantitative estimate for the sampling Kantorovich series in terms of the modulus of continuity in Orlicz spaces}, Constr. Math. Anal., \textbf{2}(1) (2019) 8-14.

\bibitem{CV2}
D. Costarelli and G. Vinti, \textit{Approximation properties of the sampling Kantorovich operators: regularization, saturation, inverse results and Favard classes in $L^{p}$-spaces}, J. Fourier Anal. Appl., \textbf{28} (2022) 49.

\bibitem{Cruz}
D. V. Cruz-Uribe and A. Fiorenza, \textit{Variable Lebesgue spaces, Applied and Numerical Harmonic Analysis}, (Birkh\"auser/Springer, Heidelberg, 2013).

\bibitem{Vore}
R. A. DeVore and G. G. Lorentz, \textit{Constructive Approximation}, Vol. 303 (Springer-Verlag, Berlin, 1993).

\bibitem{Edmund}
D. E. Edmunds and Yu. Netrusov, \textit{Entropy numbers of embeddings of Sobolev spaces in Zygmund spaces}, Studia Math., \textbf{128}(1) (1998) 71-102.


\bibitem{Gal}
S. G. Gal, \textit{A possibilistic approach of the max-product Bernstein kind operators}, Results Math., \textbf{65} (2014) 453-462.  

\bibitem{Gungor}
S. Y. G\"ung\"or and N. Ispir, \textit{Approximation by Bernstein-Chlodowsky operators of max-product kind}, Math. Commun., \textbf{23} (2018) 205-225.

\bibitem{Higgins}
J. R. Higgins, \textit{Five short stories about the cardinal series}, Bull. Amer. Math. Soc., \textbf{12} (1985) 45-89.

\bibitem{Higgins2}
J. R. Higgins, \textit{Sampling Theory in Fourier and Signal Analysis: Foundations}, Oxford Univ. Press, Oxford (1996).

\bibitem{Holhos}
A. Holhos, \textit{Weighted approximation of functions by Meyer-Konig and Zeller operators of max-product type}, Numer. Funct. Anal. Optim., \textbf{39}(6) (2018) 689-703. 

\bibitem{Kadak}
U. Kadak, \textit{Max-product type multivariate sampling operators and applications to image processing}, Chaos, Solitons and Fractals: The Interdisciplinary Journal of Nonlinear Science, and Nonequilibrium and Complex Phenomena, \textbf{157} (2022), Article 111914. 

\bibitem{Kivinukk}
A. Kivinukk and G. Tamberg, \textit{Interpolating generalized Shannon sampling operators, their norms and approximation properties}, Sampl. Theory Signal Image Process., \textbf{8} (2009) 77-95. 

\bibitem{Kivinukk2}
A. Kivinukk and G. Tamberg, \textit{On window methods in generalized Shannon sampling operators}, Appl. Numer. Harmonic Anal., \textbf{20} (2014) 63-85.  

\bibitem{Marvasti}
F. Marvasti, \textit{Nonuniform Sampling: Theory and Practice} (Springer, New York, 2001). 

\bibitem{Musielak}
J. Musielak and W. Orlicz, \textit{On modular spaces}, Studia Math, \textbf{28} (1959) 49-65.

\bibitem{Musielak2}
J. Musielak, \textit{Orlicz Spaces and Modular Spaces}, Lecture Notes in Mathematics, \textbf{1034} (Springer-Verlag, Berlin, 1983).


\bibitem{Orlova}
O. Orlova and G. Tamberg, \textit{On approximation properties of generalized Kantorovich-type sampling operators}, J. Approx. Theory, \textbf{201} (2016) 73-86.

\bibitem{Peetre}
J. Peetre, \textit{A new approach in interpolation spaces}, Studia Math., \textbf{34} (1970) 23-42. 

\bibitem{Rao}
M. Rao and Z. Ren, \textit{Applications of Orlicz spaces}, Monographs and Textbooks in Pure and Applied Mathematics, \textbf{250}, Marcel Dekker Inc., New York, 2002. 

\bibitem{Ries}
S. Ries and R. L. Stens, \textit{Approximation by generalized sampling series}, Constructive Theory of Functions'84, Sofia (1984) 746-756. 

\bibitem{Tarsi}
C. Tarsi, \textit{Adams'inequality and limiting Sobolev embeddings into Zygmund spaces}, Potential Anal., \textbf{37}(4) (2012) 353-385. 

\bibitem{Yang}
X. P. Yang, X. G. Zhou, B. Y. Cao, \textit{Single-variable term semi-latticized fuzzy relation geometric programming with max-product operator}, Inf. Sci., \textbf{325} (2015) 271-287.

\end{thebibliography}
\end{document}